\documentclass[12pt,leqno]{amsart}

\usepackage[utf8]{inputenc}   
\usepackage[T1]{fontenc}      
\usepackage{geometry}         
\usepackage[francais]{babel}  

\usepackage{amssymb}
\usepackage{enumerate}

\voffset-1cm \markleft{\rm Poisson} 
\newtheorem{theo}{Th\'{e}or\`{e}me}
\newtheorem{defi}{D\'{e}finition}
\newtheorem{propo}{Proposition}
\newtheorem{coro}{Corollaire}

\newcommand{\K}{\mathbb K}
\newcommand{\C}{\mathbb C}

\newcommand{\R}{\mathbb R}

\newcommand{\dd}{\noindent{\it D\'emonstration. }}

\newcommand{\ds}{\displaystyle}

\newcommand{\dm}{\noindent {\it D\'emonstation. }}
\title{Déterminant d'une matrice carrée}
\author{Elisabeth Remm}
\date{}
\address{Universit\'e de Haute-Alsace, IRIMAS UR 7499, F-68100 Mulhouse, France.}
\email{elisabeth.remm@uha.fr}

\begin{document}

\maketitle

Si $a$ désigne une matrice carrée à coefficients réels ou complexes, le calcul de son déterminant est initié dès les premières années du lycée, en commençant par les matrices d'ordre $2$ et le produit en croix, puis pour les matrices d'ordre $3$. Dès les premières années universitaires, la théorie du déterminant est mise en place et on expose la formule générale permettant un calcul long mais systématique que l'on rappelle ici:

Si $a=(a_{i,j})$ est une matrice carrée d'ordre $n$ à coefficients dans $\R$ ou $\C$ alors son déterminant $\det a$ est le scalaire
$$\det a=\ds \sum_{i=1}^n (-1)^{i+j}a_{i,j}M_{i,j}$$
où $M_{i,j}$ est le mineur de $a_{i,j}$, c'est-à-dire de déterminant de la matrice  d'ordre $n-1$ obtenue à partir de la matrice $a$ en lui enlevant la ligne $i$ et la colonne $j$.

Rappelons également que si $\det a \neq 0$, alors $a$ est inversible et son inverse est donnée par la formule
$$a^{-1}=\frac{1}{\det a} \-^tA$$
où $A$ est la matrice $A=((-1)^{i+j}M_{i,j})$ et où $\-^tA$ désigne la transposée de la matrice $A$.  
Nous pouvons reécrire cette formule, que $a$ soit inversible ou non, sous la forme
$$a \cdot \-^tA= \det(a) I_n$$
où $I_n$ désigne la matrice identité d'ordre $n$.

Dans ce petit mémoire, nous allons voir d'autres formules permettant ce calcul du déterminant, dont l'origine remonte à la résolution d'un exercice issu du livre

http:://livres-mathematiques.fr/algebre-multilineaire

\section{Quelques rappels et matrice des cofacteurs}

Si $a=(a_{i,j})$ est une matrice carrée d'ordre $n$, on appelle mineur du couple $(i,j)$ le déterminant de la matrice où on a barré la 
$i$-ème ligne et la $j$-ème colonne. Si ce mineur est noté $M_{i,j}$, le cofacteur du couple $(i,j)$  est $A_{i,j}=(-1)^{i+j} M_{i,j}.$ Soit $A$ la matrice $A=((-1)^{i+j}M_{i,j})$. Nous avons vu que 
$$a \cdot  \-^tA=\det a . Id_n.$$
Si $a$ est inversible, on retrouve la formule donnant l'inverse de $a$. Si $a$ est singulière, on a alors l'identité
 $$a.\-^tA=0.$$
Dans ce cas, si $A \neq 0$, ce qui est équivalent à dire que le rang de $a$ est $n-1$, alors les lignes de $A$ engendrent le noyau de $a$ qui est de dimension $1$ et toutes les lignes de $A$ sont proportionnelles. Si $A=0$, le rang de $a$ est inférieur ou égal à $n-2$. Dans ce cas, nous sommes conduits à étudier la matrice des "cofacteurs" d'ordre $n-2$. 
 
 \noindent Exemple. Soit la matrice 
 $$a
=\begin{pmatrix}
     1 &-1 & 1 &0    \\
   0   & 1 & 1 & 1\\
   0 &0 & 0 & 0\\ 
   1 & 1 & 3 & 2
\end{pmatrix}
$$
On a 
$$A=0$$
et le rang de $a$ est au plus égal à $2$.  Nous allons voir dans le paragraphe suivant, le rôle des mineurs d'ordre $2$.

\section{Déterminant d'une matrice d'ordre $4$}

Soit à calculer le déterminant de la matrice 
$$a=\begin{pmatrix}
    a_{1,1}  &  a_{1,2}  &a_{1,3}  &a_{1,4}    \\
      a_{2,1}  &  a_{2,2}  &a_{2,3}  &a_{2,4}    \\  
        a_{3,1}  &  a_{3,2}  &a_{3,3}  &a_{3,4}    \\    
        a_{4,1}  &  a_{4,2}  &a_{4,3}  &a_{4,4}    \\
\end{pmatrix}
$$
Notons par $a_{(i_1i_2,j_1,_2)}$, $i_1 <i_2$ et $j_1 <j_2$ la matrice obtenue à partir de la matrice $a$ en ne conservant que les lignes $i_1$ et $i_2$ et les colonnes $j_1$ et $j_2$. Nous noterons également par  $a_{(\overline{i_1i_2,j_1,j_2})}$, $i_1<i_2, \ j_1<j_2$ la matrice obtenue à partir de $a$ en enlevant les lignes $i_1$ et $i_2$ et les colonnes $j_1$ et $j_2$.  Par exemple
$$a_{(12,12)}=\begin{pmatrix}
      a_{1,1} &    a_{1,2} \\
      a_{2,1} &    a_{2,2} \\
\end{pmatrix}\ \ a_{(\overline{12,12})}=\begin{pmatrix}
     a_{3,3} &    a_{3,4} \\
      a_{4,3} &    a_{4,4} \\
\end{pmatrix}, \  \ \ a_{(\overline{34,34})}=\begin{pmatrix}
     a_{1,1} &    a_{1,2} \\
      a_{2,1} &    a_{2,2} \\
\end{pmatrix}.$$
Bien entendu $a_{(\overline{34,34})}=a_{(12,12)}$. 
Posons
$$M_{i_1i_2}^{j_1j_2}=\det a_{(i_1i_2,j_1,j_2)}, \ \ M_{(\overline{i_1i_2}}^{\overline{j_1j_2}}=\det a_{(\overline{(i_1i_2,j_1,j_2})}$$

On dira que $M_{(\overline{i_1i_2}}^{\overline{j_1j_2}}$ est 
un mineur  d'ordre $2$ de la matrice $a$. 
On a alors
\begin{theo}
$$\det a= \ds \sum_{(j_1,j_2)=(1,2)}^{(3,4)} (-1)^{i_1+i_2+j_1+j_2}M_{i_1i_2}^{j_1j_2}M_{(\overline{i_1i_2}}^{\overline{j_1j_2}}.$$
\end{theo}
Dans cette formule on suppose que les couples $(j_1,j_2)$ prennent successivement les valeurs 
$$(1,2),(1,3),(1,4),(2,3),(2,4),(3,4).$$
\dm Comme le déterminant est invariant,au signe près, de l'ordre des lignes, nous pouvons supposer, pour simplifier l'écriture, que $(i_1,i_2)=(1,2).$ Nous devons donc vérifier que 
$$\det a= M_{12}^{12}M_{34}^{34}-M_{12}^{13}M_{34}^{24}+M_{12}^{14}M_{34}^{23}+M_{12}^{23}M_{34}^{14}-M_{12}^{24}M_{34}^{13}+M_{12}^{34}M_{34}^{12}.$$
Ce calcul ne présente aucune difficulté. En effet
$$
\begin{array}{ll}
\det a = & a_{1,4} a_{2,3} a_{3,2} a_{4,1} - a_{1,3 } a_{2,4 } a_{3,2} a_{4,1} - a_{1,4 } a_{2,2 } a_{3,3 } a_{4,1} + 
 a_{1,2 } a_{2,4 } a_{3,3 } a_{4,1} \\
 &+ a_{1,3 } a_{2,2 } a_{3,4 } a_{4,1} - a_{1,2 } a_{2,3 } a_{3,4 } a_{4,1} - 
 a_{1,4 } a_{2,3 } a_{3,1 } a_{4,2} + a_{1,3 } a_{2,4 } a_{3,1 } a_{4,2 }\\
 &+ a_{1,4 } a_{2,1 } a_{3,3 } a_{4,2} - 
 a_{1,1 } a_{2,4 } a_{3,3 } a_{4,2}- a_{1,3 } a_{2,1 } a_{3,4 } a_{4,2 }+ a_{1,1 } a_{2,3 } a_{3,4 } a_{4,2 }\\
 &+ 
 a_{1,4 } a_{2,2 } a_{3,1 } a_{4,3} - a_{1,2 } a_{2,4 } a_{3,1 } a_{4,3} -  a_{1,4 } a_{2,1 } a_{3,2} a_{4,3} + 
 a_{1,1 } a_{2,4 } a_{3,2 } a_{4,3} \\
 &+ a_{1,2 } a_{2,1 } a_{3,4 } a_{4,3 }- a_{1,1 } a_{2,2 } a_{3,4 } a_{4,3} - 
 a_{1,3 } a_{2,2 } a_{3,1 } a_{4,4} + a_{1,2 } a_{2,3 } a_{3,1 } a_{4,4} \\
 &+ a_{1,3 } a_{2,1 } a_{3,2 } a_{4,4} - 
 a_{1,1 } a_{2,3 } a_{3,2 } a_{4,4} - a_{1,2 } a_{2,1 } a_{3,3 } a_{4,4} + a_{1,1 } a_{2,2 } a_{3,3 } a_{4,4}.
 \end{array}
 $$
 ce qui correspond au développement du second membre. 
 
 Nous pouvons associer à la matrice $a$, la matrice des cofacteurs d'ordre $2$. Plus précisément, notons par $m_2(a)$ la matrice
 $$m_2(a)=\begin{pmatrix}
     M_{12}^{12} &   M_{12}^{13} &M_{12}^{14} &M_{12}^{23} &M_{12}^{24} &M_{12}^{34}  \\
        M_{13}^{12} &   M_{13}^{13} &M_{13}^{14} &M_{13}^{23} &M_{13}^{24} &M_{13}^{34}  \\
           M_{14}^{12} &   M_{14}^{13} &M_{14}^{14} &M_{14}^{23} &M_{14}^{24} &M_{14}^{34}  \\
   M_{23}^{12} &   M_{23}^{13} &M_{23}^{14} &M_{23}^{23} &M_{23}^{24} &M_{23}^{34}  \\
      M_{24}^{12} &   M_{24}^{13} &M_{24}^{14} &M_{24}^{23} &M_{24}^{24} &M_{24}^{34}  \\
         M_{34}^{12} &   M_{34}^{13} &M_{34}^{14} &M_{34}^{23} &M_{34}^{24} &M_{34}^{34}  \\
\end{pmatrix}
$$
et soit la matrice $\widetilde{m_2(a)}$ définie par
$$\widetilde{m_2(a)}=\begin{pmatrix}
     M_{34}^{34} &   -M_{24}^{34} &M_{23}^{34} &M_{14}^{34} &-M_{13}^{34} &M_{12}^{34}  \\
           -M_{34}^{24} &   M_{24}^{24} &-M_{23}^{24} &-M_{14}^{24} &M_{13}^{24} &-M_{12}^{24}  \\
               M_{34}^{23} &   -M_{24}^{23} &M_{23}^{23} &M_{14}^{23} &M_{13}^{23} &M_{12}^{23}  \\
               M_{34}^{14} &   -M_{24}^{14} &M_{23}^{14} &M_{14}^{14} &M_{13}^{14} &M_{12}^{14}  \\
                    -M_{34}^{13} &   M_{24}^{13} &-M_{23}^{13} &-M_{14}^{13} &-M_{13}^{13} &-M_{12}^{13}  \\
                         M_{34}^{12} &   -M_{24}^{12} &M_{23}^{12} &M_{14}^{12} &M_{13}^{12} &M_{12}^{12}  \\
\end{pmatrix}
$$
Cette matrice correspond au produit
$$I_{6,2}\cdot J \cdot m_2(a) \cdot J \cdot I_{6,2}$$
où $$J=\begin{pmatrix}
   0 &0& 0&0&0&1\\
   0 &0& 0&0&1& 0\\   
    0 &0& 0&1&0& 0\\ 
    0 &0& 1&0&0& 0\\ 
    0 &1& 0&0&0& 0\\ 
    1 &0& 0&0&0& 0\\    
\end{pmatrix}$$
et 
$I_{6,2}$ la matrice diagonale $Diag(1,-1,1,1,-1,1)$.
On vérifie la relation
\begin{equation}
\label{m2}
m_2(a)\cdot \widetilde{m_2(a)}= \det a . Id_6
\end{equation}
En particulier, si $a$ est une matrice singulière, alors $m_2(a)$ est aussi singulière 
$$m_2(a)\cdot \widetilde{m_2(a)}= 0$$
sinon l'inverse de la matrice $m_2(a)$ est $\frac{1}{\det a}\widetilde{m_2(a)}.$.

\noindent{\bf Exemple.} 
Revenons à l'exemple exposé au paragraphe précédent. On a 
$$m_2(a)=\begin{pmatrix}
    1 &   1 &1 &-2 &-1 &1  \\
       0 &   0&0&0&0&0 \\
        2 &   2 &2&-4 &-2&2 \\
              0 &   0&0&0&0&0 \\
                                 -1 &   -1 &-1 &2 &1&-1  \\
                       0 &   0&0&0&0&0 \\
                       \end{pmatrix}
$$
et cette matrice est de rang $1$ et
$$\widetilde{m_2(a)}=
\begin{pmatrix}
     0 & 1 & 0 & 2 & 0 & 1   \\
       0 & 1 & 0 & 2 & 0 & 1   \\
         0 & -2 & 0 & -4& 0 & -2 \\
           0 & 1 & 0 & 2 & 0 & 1   \\
             0 & -1 & 0 & -2 & 0 & -1   \\
               0 & 1 & 0 & 2 & 0 & 1   \\
\end{pmatrix}
$$
et $m_2(a) \cdot \widetilde{m_2(a)}=0.$
On vérifie que l'on a l'identité matricielle
$$
a\cdot
\begin{pmatrix}
     M_{12}^{23} &M_{12}^{24} &M_{12}^{34} &0  \\
    -M_{12}^{13} &-M_{12}^{14}&0 &M_{12}^{34} \\
    M_{12}^{12} &   0&-M_{12}^{14} &-M_{12}^{24}\\
    0 & M_{12}^{12} &   M_{12}^{13} &M_{12}^{23}  \\
\end{pmatrix}=
$$
$$\begin{pmatrix}
     1 &-1 & 1 &0    \\
   0   & 1 & 1 & 1\\
   0 &0 & 0 & 0\\ 
   1 & 1 & 3 & 2
\end{pmatrix}
\cdot 
\begin{pmatrix}
    -2 & -1 & 1 & 0\\
    -1 & -1 & 0 & 1 \\
    1 & 0 & -1 & 1\\
    0 & 1 & 1 & -2
    \end{pmatrix}=0
$$
Ainsi les vecteurs colonnes de la matrice $$\begin{pmatrix}
     M_{12}^{23} &M_{12}^{24} &M_{12}^{34} &0  \\
    -M_{12}^{13} &-M_{12}^{14}&0 &M_{12}^{34} \\
    M_{12}^{12} &   0&-M_{12}^{14} &-M_{12}^{24}\\
    0 & M_{12}^{12} &   M_{12}^{13} &M_{12}^{23}  \\
\end{pmatrix},$$
que nous appellerons la matrice des cofacteurs d'ordre $2$ basée sur $(i_1,i_2)=(1,2)$, 
engendre le noyau de $a$. Nous avons donc montré que si $a$ est singulière, la matrice des cofacteurs détermine le noyau de $a$ si le corang de $a$ est $1$ sinon la matrice des cofacteurs d'ordre $2$ détermine ce noyau si le corang est égal à $2$.

\medskip

Ceci se généralise ainsi
\begin{propo}
Si $a$ est une matrice d'ordre $4$ de rang $2$, alors pour tout couple $(i_1,i_2)$ avec $i_1<i_2$, on a 
$$
a\cdot
\begin{pmatrix}
         M_{i_1i_2}^{23} &M_{i_1i_2}^{24} &M_{i_1i_2}^{34} &0  \\
    -M_{i_1i_2}^{13} &-M_{i_1i_2}^{14}&0 &M_{i_1i_2}^{34} \\
    M_{i_1i_2}^{12} &   0&-M_{i_1i_2}^{14} &-M_{i_1i_2}^{24}\\
    0 & M_{i_1i_2}^{12} &   M_{i_1i_2}^{13} &M_{i_1i_2}^{23}  \\\end{pmatrix}=0.
$$
\end{propo}
En effet, les coefficients non nuls correspondent aux mineurs d'ordre $3$ qui sont nuls par hypothèse.

\medskip

Ceci étant, la relation (\ref{m2}) implique
$$\det(m_2(a))^2=\det (a)^6.$$
On en déduit, après vérification des signes
\begin{propo}
$$\det (m_2(a))=(\det(a))^3.$$
\end{propo}
Il est facile de visualiser cette proprité lorsque $a$ est une matrice triangulaire. 

\medskip

Supposons $a \in SL(4)$, c'est-à-dire $\det a =1$. D'après la proposition précédente $\det(m_2(a)=1$ et donc
$m_2(a) \in SL(6).$
On a donc 
$$a \in SL(4) \Longrightarrow m_2(a) \in SL(6).$$
Ceci conduit à examiner l'application $\varphi_2 : \mathcal{M}(4,\K) \rightarrow \mathcal{M}(6,\K)$ définie par
$$\varphi_2 (a) =m_2(a).$$
D'après la proposition précédente, cette application induit une application, que l'on note toujours $\varphi_2$:
$$\varphi_2: GL(4,\K) \rightarrow GL(6,\K)$$
avec le notations, $\mathcal{M}(n,\K)$ désigne l'espace vectoriel des matrices carrées d'ordre $n$ et $GL(n,\K)$ le groupe des matrices inversibles d'ordre $n$. Rappelons que  $A_{ik}^{js}$ désigne le mineur d'ordre $2$ issu de la matrice $a$ et posé sur les lignes $i$ et $ k$ et sur les colonnes $j$ et $s$ , alors si $c=ab$ est le produit des matrices $a$ et $b$, et si $M(a)^{i,j}_{k,l}$ désigne le mineur du tyê $(i,j),(k,l)$ associé à la matrice $a$, alors
$$
\begin{array}{ll}
   M(c)_{ik}^{js}=  & M(a)_{ik}^{12}M(b)_{12}^{js}+M(a)_{ik}^{13}M(b)_{13}^{js}+M(a)_{ik}^{14}M(b)_{14}^{js}  
 + M(a)_{ik}^{23}M(b)_{23}^{js}\\&+M(a)_{ik}^{24}M(b)_{24}^{js}
+M(a)_{ik}^{34}M(b)_{34}^{js}
 \end{array}
$$
Cette formule se déduit directement de l'écriture
$$ab=(<L_ia,C_jb>)$$
où $L_i(a)$ désigne le vecteur ligne $i$ de $a$ et $C_j(b)$ le vecteur colonne $j$ de $b$ ce qui implique
$$C_{ik}^{js}=<L_ia,C_jb><L_ka,C_sb>-<L_ia,C_sb><L_ka,C_jb>.$$
Nous obtenons ainsi les coefficients de la matrice $\varphi_2(ab)=m_2(ab).$ La matrice $m_2(a)m_2(b)$ est donnée par
$$m_2(a)m_2(b)=( \sum_{kl \in I}M(a)_{23}^{kl}M(b)_{kl}^{13})$$

{\tiny
$$ \begin{pmatrix}
 \ds   \sum_{ij \in I}
    M(a)_{12}^{ij}M(b)_{ij}^{12}  &    \ds  \sum_{ij \in I}M(a)_{12}^{ij}M(b)_{ij}^{13}  &  \ds \sum_{ij \in I}M(a)_{12}^{ij}M(b)_{ij}^{14}  &  \ds \sum_{ij \in I}M(a)_{12}^{ij}M(b)_{ij}^{23}  & \ds  \sum_{ij \in I}M(a)_{12}^{ij}M(b)_{ij}^{24}  &  \ds \sum_{ij \in I}M(a)_{12}^{ij}M(b)_{ij}^{34}   \\
     \ds  \sum_{ij \in I} M(a)_{13}^{ij}M(b)_{ij}^{12}  &    \ds  \sum_{ij \in I}M(a)_{13}^{ij}M(b)_{ij}^{13}  & \ds  \sum_{ij \in I}M(a)_{13}^{ij}M(b)_{ij}^{14}  &  \ds  \ds \sum_{ij \in I}M(a)_{13}^{ij}M(b)_{ij}^{23}  &  \ds \sum_{ij \in I}M(a)_{13}^{ij}M(b)_{ij}^{24}  & \ds  \sum_{ij \in I}M(a)_{13}^{ij}M(b)_{ij}^{34}   \\ 
        \ds \sum_{ij \in I}  M(a)_{14}^{ij}M(b)_{ij}^{12}  &    \ds \sum_{ij \in I}M(a)_{14}^{ij}M(b)_{ij}^{13}  & \ds \sum_{ij \in I}M(a)_{14}^{ij}M(b)_{ij}^{14}  & \ds \sum_{ij \in I}M(a)_{14}^{ij}M(b)_{ij}^{23}  & \ds \sum_{ij \in I}M(a)_{14}^{ij}M(b)_{ij}^{24}  & \ds \sum_{ij \in I}M(a)_{14}^{ij}M(b)_{ij}^{34}   \\
             \ds \sum_{ij \in I} M(a)_{23}^{ij}M(b)_{ij}^{12}  &    \ds \sum_{ij \in I}M(a)_{23}^{ij}M(b)_{ij}^{13}  & \ds \sum_{ij \in I}M(a)_{23}^{ij}M(b)_{ij}^{14}  & \ds \sum_{ij \in I}M(a)_{23}^{ij}M(b)_{ij}^{23}  & \ds \sum_{ij \in I}M(a)_{23}^{ij}M(b)_{ij}^{24}  & \ds \sum_{ij \in I}M(a)_{23}^{ij}M(b)_{ij}^{34}   \\
                 \ds \sum_{ij \in I} M(a)_{24}^{ij}M(b)_{ij}^{12}  &    \ds \sum_{ij \in I}M(a)_{24}^{ij}M(b)_{ij}^{13}  & \ds \sum_{ij \in I}M(a)_{24}^{ij}M(b)_{ij}^{14}  & \ds \sum_{ij \in I}M(a)_{24}^{ij}M(b)_{ij}^{23}  & \ds \sum_{ij \in I}M(a)_{24}^{ij}M(b)_{ij}^{24}  & \ds \sum_{ij \in I}M(a)_{24}^{ij}M(b)_{ij}^{34}   \\
                     \ds \sum_{ij \in I} M(a)_{34}^{ij}M(b)_{ij}^{12}  &    \ds \sum_{ij \in I}M(a)_{34}^{ij}M(b)_{ij}^{13}  & \ds \sum_{ij \in I}M(a)_{34}^{ij}M(b)_{ij}^{14}  & \ds \sum_{ij \in I}M(a)_{34}^{ij}M(b)_{ij}^{23}  & \ds \sum_{ij \in I}M(a)_{34}^{ij}M(b)_{ij}^{24}  & \ds \sum_{ij \in I}M(a)_{34}^{ij}M(b)_{ij}^{34}   \\
\end{pmatrix}
$$
}
où $I=\{12,13,14,23,24,34\}$. Chacun de ces coefficients correspond bien au coefficient $M(b)_{ik}^{jl}$. On en déduit
$$m_2(a)m_2(b)=m_2(ab).$$
\begin{propo}
L'application
$$\varphi_2:GL(4,\K) \rightarrow GL(6,\K)$$
définie par $$\varphi_2 (a)=m_2(a)$$
est un homomorphisme de groupes.
\end{propo}
Bien entendu cet homomorphisme n'est pas surjectif, ne serait-ce que par les dimensions (des variétés différentielles) des ensembles de départ et d'arrivée. Par exemple, si 
$$b=Diag[\mu_1,\mu_2,\mu_3,\mu_4,\mu_5,\mu_6]$$
 avec le produit de ces $\mu_i$ non nul, alors il existe $a \in GL(4,\K)$ telle que $b=m_2(a)$ si et seulement si $\mu_1\mu_6=\mu_2\mu_5=\mu_3\mu_4$. Dans ce cas, la matrice $a$ s'écrit
$$a=\begin{pmatrix}
      a_{11}& 0 & 0 & 0   \\
      0&  \mu_1/a_{11} & 0 & 0 \\
      0 & 0&  \mu_2/a_{11} & 0  \\
       0& 0 & 0 &  \mu_3/a_{11} \\
\end{pmatrix}
$$
avec $a_{11} \neq 0$ et $a_{11} ^2=\frac{\mu_1\mu_2}{\mu_4}=\frac{\mu_1\mu_3}{\mu_5}=\frac{\mu_2\mu_3}{\mu_6}.$ Par contre l'application $\varphi_2$ est injective. Pour cela considérons cette application comme une application différentiable entre les deux variétés $GL(4,\K)$ et $GL(6,\K)$. Sa matrice jacobienne en l'élément neutre de $GL(4,\K)$ est de rang $16$ et le théorème du rang nous assure que $\varphi_2$ est injective.
\begin{coro}
L'application $\varphi_2$ est un homomorphisme injectif.
\end{coro}
On en déduit la suite exacte
$$0 \rightarrow GL(4,\K) \rightarrow GL(6,\K).$$
Nous avons vu précédemment que si $a$ est de déterminant $1$, il en est de même de $m_2(a)$. Ceci se déduit aussi du corallaire ci-dessus. Il en est également pour les matrices orthogonales
\begin{propo}
Soit $a \in GL(4,\K)$. Alors \begin{enumerate}
  \item Si $a \in SL(4,\K)$ alors $m_2(a) \in SL(6,\K)$.
  \item Si $a \in SO(4,\K)$ alors $m_2(a) \in S0(6,\K)$.
\end{enumerate}
\end{propo}

\medskip

\noindent{\bf Remarque.} La formule classique du calcul du déterminant de $a$ utilise la matrice des mineurs d'ordre $3$. 
En respectant les notations ci-dessus, nous allons noter par $M_{ijk}^{lmn}$ le déterminant de la matrice $a_{(ijk,lmn)}$ obtenue à partir de la matrice $a$ en ne conservant que les lignes $i,j,k$ et les colonnes $l,m,n$ avec $i<j<k$ et $l<m<n$ (avec les notations de l'introduction, nous avons par exemple $M_{123}^{123}=M_{44}$). Notons par $m_3(a)$ la matrices
$$
m_3(a)=\begin{pmatrix}
     M_{123}^{123} &   M_{124}^{124} &M_{123}^{134} &M_{123}^{234} \\
        M_{124}^{123} &   M_{124}^{124} &M_{124}^{134} &M_{124}^{234}  \\
           M_{134}^{123} &   M_{134}^{124} &M_{134}^{134} &M_{134}^{234}  \\
   M_{234}^{123} &   M_{234}^{124} &M_{234}^{134} &M_{234}^{234} 
\end{pmatrix}
$$
et par $\widetilde{m_3(a)}$ la matrice:
$$\widetilde{m_3(a)}=I_{4,2}.J.\-^tm_3(a). J.I_{4,2}$$
où $I_{4,2}$ est la matrice diagonale $(1,-1,1,-1)$ c'est-à-dire
$$\widetilde{m_3(a)}=
\begin{pmatrix}
     M_{234}^{234} &   -M_{134}^{234} &M_{124}^{234} &-M_{123}^{234} \\
        -M_{234}^{134} &   M_{134}^{134} &-M_{124}^{134} &M_{123}^{134}  \\
           M_{234}^{124} &   -M_{134}^{124} &M_{124}^{124} &-M_{123}^{124}  \\
   -M_{234}^{123} &   M_{134}^{123} &-M_{124}^{123} &M_{123}^{123} 
\end{pmatrix}
$$

on a la formule classique
$$a \cdot \widetilde{m_3(a)}=\det (a) I_4.$$
On en déduit
$$\det m_3(a)=(\det (a))^3.$$
En particulier si $a \in SL(4)$ alors $m_3(a) \in SL(4).$

\begin{propo}
L'application $\varphi_2 : GL(4) \rightarrow GL(4)$ définie par $\varphi_2 (a)=m_3(a)$ est un homomorphisme de groupe.
\end{propo}

En effet, si $a$ est non singulière, alors $a^{-1}=\frac{1}{\det(a)}\widetilde{m_3(a)}$. Ansi
$$\widetilde{m_3(a)}=\det(a) a^{-1}$$ 
ce qui implique
$$m_3(a)=(\det(a))I_{4,2}.J.\-^ta^{-1}. J.I_{4,2}$$
Ainsi
$$m_3(ab)=\det(ab)I_{4,2}.J.\-^ta^{-1}. J.I_{4,2}.I_{4,2}.J.\-^tb^{-1}. J.I_{4,2}=\det(a)\det(b)\-^ta^{-1}\-^tb^{-1}=m_3(a)m_3(b)$$
pour toutes matrices $a,b \in GL(4).$

\begin{coro}
Soit $a \in GL(4,\K)$. Alors \begin{enumerate}
  \item Si $a \in SL(4,\K)$ alors $m_3(a) \in SL(4,\K)$.
  \item Si $a \in SO(4,\K)$ alors $m_3(a)\in S0(4,\K)$.
\end{enumerate}
\end{coro}
Examinons le morphisme $\varphi_2: SO(4) \rightarrow SO(4)$. On a $\varphi_2(a)=m_3(a)$. Or comme $a \in SO(4)$, $a^{-1}=\-^ta.$ On en déduit
$$m_3(a)=(\det(a))I_{4,2}.J.\-^ta^{-1}. J.I_{4,2}=I_{4,2}.J.a. J.I_{4,2}$$
et comme $m_3(a) \in SO(4)$, 
$$m_3(m_3(a))=I_{4,2}.J.m_3(a). J.I_{4,2}=I_{4,2}.J.I_{4,2}.J.a. J.I_{4,2}. J.I_{4,2}=a$$
 et donc
$$m_3(m_3(a))=a.$$
L'homomorphisme $\varphi_2$ restreint à $SO(4)$ vérifie
$$\varphi_2 \circ \varphi_2 =Id.$$
 Dans le cas général, on a 
\begin{propo}
 Soit $a \in GL(4,\K)$. Alors
 $$\varphi_2^2(a)=m_3(m_3(a))=(\det (a))^2 a.$$
 En particulier, $\varphi_2$ est une symétrie sur $SL(4,\K)$ et sur tous ses sous-groupes .
 \end{propo}
 \dd De la formule donnant $a^{-1}$, nous avons déduit
$$m_3(a)=(\det(a))I_{4,2}.J.\-^ta^{-1}. J.I_{4,2}.$$
Ceci implique
$$
\begin{array}{ll}
m_3(m_3(a))&=(\det (m_3(a))((I_{4,2}.J. ^t(m_3(a)^{-1}) . J.I_{4,2})\\
&=\frac{\det m_3(a)}{\det(a)}a.
\end{array}
$$
Comme $\det(m_3(a))=(\det(a))^3$ et 
$$m_3(m_3(a))=(\det(a))^2) a.$$
D'où la proposition.

\section{Cas général}
Soit $a$ une matrice carrée d'ordre $n$. Notons par $m_p(a)$ la matrice des mineurs d'ordre $p$
$$m_p(a)=(M_{i_1,\cdots,i_p}^{j_1,\cdots ,j_p})$$
où $M_{i_1,\cdots,i_p}^{j_1,\cdots ,j_p}$ est le déterminant de la matrice d'ordre $p$ obtenue en ne conservant que les lignes $i_1,\cdots,i_p$ et les colonnes $j_1,\cdots,j_p$ avec $i_1 < \cdots <i_p, \ j_1 < \cdots < j_p$. Pour ranger les éléments de cette matrice on utilise l'ordre naturel lexicographique. Par exemple, si $n=4$ et $p=3$, on aura
$$m_3(a)=\begin{pmatrix}
    M_{1,2,3}^{1,2,3}  &  M_{1,2,3}^{1,2,4} & M_{1,2,3}^{1,3,4} & M_{1,2,3}^{2,3,4}  \\
     M_{1,2,4}^{1,2,3}  &  M_{1,2,4}^{1,2,4} & M_{1,2,4}^{1,3,4} & M_{1,2,4}^{2,3,4}  \\
       M_{1,3,4}^{1,2,3}  &  M_{1,3,4}^{1,2,4} & M_{1,3,4}^{1,3,4} & M_{1,3,4}^{2,3,4}  \\
         M_{2,3,4}^{1,2,3}  &  M_{2,3,4}^{1,2,4} & M_{2,3,4}^{1,3,4} & M_{2,3,4}^{2,3,4}  \\ 
\end{pmatrix}
$$
La matrice $m_p(a)$ est d'ordre $\binom{n}{p}$. A cette matrice on associe la matrice $\overline{m_p(a)}$ dont les éléments sont les mineurs $A_{\overline{i_1,\cdots,i_p}}^{\overline{j_1,\cdots ,j_p}}$ où $\overline{i_1,\cdots,i_p}$ est le complémentaire dans $\left[\left[1,\cdots,n\right]\right]$ de $\{i_1,\cdots,i_p\}$. Enfin, nous noterons par $\widetilde{m_p(a)}$ la matrice
$$\widetilde{m_p(a)}=Diag(1,-1,1,\cdots) \overline{m_p(a)} Diag(1,-1,1,\cdots)$$
où $Diag(u_1,u_2,\cdots)$ désigne la matrice diagonale dont les éléments sur la diagonales sont dans l'ordre $u_1,u_2,\cdots.$

\begin{propo}
Soit $a$ une matrice carrée d'ordre $n$. Alors
$$m_p(a) \cdot \widetilde{m_p(a)}= \det (a) Id_{n_p}$$
où $n_p=\binom{n}{p}.$
\end{propo}
\subsection{Mineurs d'ordre $n-1$}
Soit $a$ une matrice d'ordre $n$ et $m_{n-1}$ la matrice des mineurs d'ordre $n-1$. Cette matrice est aussi d'ordre $n$ et vérifie la formule classique
$$a. \widetilde{m_{n-1}}(a)=\det (a) Id.$$
ce qui implique
$$\det(m_{n-1}(a))=(\det(a))^{n-1}.$$

Si $a$ est singulière il en est de même de $m_{n-1}(a)$ et cette dernière matrice est nulle si le rang de $a$ est inférieur ou égal à $n-2$. Si ce rang est égal à $n-1$, alors $m_{n-1}(a) \neq 0$ et ses vecteurs colonnes engendrent le noyau de $a$ qui est de dimension $1$, ainsi toutes les colonnes de $m_{n-1}$ sont proportionnelles et son rang est égal à $1$. 

Supposons $a$ non singulière, soit $\det (a) \neq 0.$
Comme nous l'avons vu au paragraphe précédent, 
$$m_{n-1}(a)=\det (a) (I_{n,2} \cdot J \cdot \-^ta^{-1} \cdot J \cdot I_{n,2})$$
où $I_{n,2}$ est la matrice diagonale d'ordre $n$ valant $Diag[1,-1,1,\cdots,(-1)^{i+1},\cdots,(-1)^{n+1}]$.  On en déduit, comme $m_{n-1}$ est aussi régulière
$$ m_{n-1}(m_{n-1}(a))     =\det (m_{n-1}(a)) (I_{n,2} \cdot J \cdot \-^tm_{n-1}(a)^{-1} \cdot J \cdot I_{n,2}).$$
Mais
$$\-^tm_{n-1}(a)^{-1}=(\det (a))^{-1} (I_{n,2} \cdot J \cdot a \cdot J \cdot I_{n,2})$$
d'où
$$\begin{array}{ll}
\medskip
 m_{n-1}(m_{n-1}(a))  &=\det (m_{n-1}(a)) (I_{n,2} \cdot J \cdot \-^tm_{n-1}(a)^{-1} \cdot J \cdot I_{n,2})\\
\medskip
    &   =\det (m_{n-1}(a)) (I_{n,2} \cdot J \cdot ( (\det (a)^{-1} (I_{n,2} \cdot J \cdot a \cdot J \cdot I_{n,2})) \cdot J \cdot I_{n,2})  \\
\medskip
& =\det (m_{n-1}(a)) (\det (a)^{-1}  a   \\
\medskip
& =(\det(a))^{n-2} a.
\end{array}
$$
Nous avons également, pour toutes matrices $a,b$ carrées d'ordre $n$
$$\begin{array}{ll}
\medskip
m_{n-1}(a)\cdot m_{n-1}(b) & =\det (a) (I_{n,2} \cdot J \cdot \-^ta^{-1} \cdot J \cdot I_{n,2}) \cdot\det (b) (I_{n,2} \cdot J \cdot \-^tb^{-1} \cdot J \cdot I_{n,2})\\
\medskip
& = \det(a)\det(b) (I_{n,2} \cdot J \cdot \-^ta^{-1}\cdot \-^tb^{-1}  \cdot J \cdot I_{n,2}\\
\medskip
&= \det(a)\det(b) (I_{n,2} \cdot J \cdot \-^t(a\cdot b)^{-1}  \cdot J \cdot I_{n,2}\\
& = m_{n-1}(a \cdot b).$$
\end{array}
$$
\begin{propo}
Soit $\varphi_{n-1}$ l'application qui à toute matrice carrée d'ordre $n$ fait correspondre la matrice carrée $m_{n-1}(a)$. Notons également par $\varphi_{n-1}$  sa restriction au groupe $GL(n,\K)$ des matrices inversibles. Alors
\begin{enumerate}
  \item $\varphi_{n-1} : GL(n,\K) \rightarrow GL(n,\K)$ est un homomorphisme de groupes.
  \item Cet homomorphisme vérifie $\varphi_{n-1} \circ \varphi_{n-1} (a) =(\det(a))^{n-1} a.$
  \item Les restrictions de cet homomorphisme aux sous-groupes $SL(n,\K)$ et $SO(n,\K)$ sont des symétries.
\end{enumerate}
\end{propo}

\subsection{Mineurs d'ordre $n-2$}

Comme précédemment, notons par $M_{i,j}^{k,l}$, $i<j,\ k<l$ le déterminant de la matrice d'ordre $2$ notée $a_{i,j}^{k,l}$ obtenue à partir de la matrice $a$ en ne conservant que les deux lignes $i$ et $j$ et les deux colonnes $k$ et $l$.  Si $\overline{i,j}$ désigne le complémentaire ordonné de ${i,j}$ dans $\left[\left[1,\cdots,n\right]\right]$, le mineur de $a_{i,j}^{k,l}$  est noté $M_{\overline{i,j}}^{\overline{k,l}}$.  C'est un déterminant d'ordre $n-2$.
\begin{propo}
$$\det(a)=\sum_{j<k} (-1)^{j+k+1} M_{1,2}^{j,k} M_{\overline{1,2}}^{\overline{j,k}}$$
\end{propo}
Comme dans le cas particulier précédent, cette identité peut s'écrire sous forme matricielle
$$m_2(a)\cdot \widetilde{m_{n-2}(a)}=\det(a) I_n$$
avec
$$\widetilde{m_{n-2}(a)}=I_{n,2}\cdot J \cdot m_{n-2}(a) \cdot J \cdot I_{n,n-2}$$ et $I_{n,2}=Diag( (-1)^{i+j+1}), 1 \leq i <j \leq n).$ En particulier, $m_2(a)$ est inversible si et seulement $a$ aussi et son inverse est $\frac{1}{\det (a)}\widetilde{m_{n-2}(a)}.$ On en déduit également
$$(\det(m_2(a))^2=(\det(a))^n$$
et donc
\begin{equation}
\label{det2}
\det(m_2(a))=(\det(a))^{n-1}.
\end{equation}
Considérons l'application
$$\varphi_2: GL(n,\K) \rightarrow GL(N,\K)$$
avec $N=\binom{n}{2}$ définie par
$$\varphi_2 (a)=m_2(a).$$
Avec les mêmes arguments que ceux développés dans le cas particulier $n=4$, on peut affirmer
\begin{propo}
L'application
$$\varphi_2: GL(n,\K) \rightarrow GL(N,\K)$$
vérifie 
$$\varphi_2(ab)=\varphi_2(a)\varphi_2(b)$$
pour toutes matrices 
$a,b \in GL(n,\K)$
\end{propo}
On en déduit que les restrictions de cet homomorphisme $SL(n,\K)$ et $SO(n,\K)$ sont aussi des homorphismes de groupes. Ceci permet de définir par exemple, la suite
$$SL(4,\K) \stackrel{\varphi_2}\rightarrow  SL(6,\K) \stackrel{\varphi_2}\rightarrow SL(10,\K) \stackrel{\varphi_2}\rightarrow \cdots $$
Par contre,nous avons
$$GL(3,\K) \stackrel{\varphi_2}\rightarrow  GL(3,\K)$$
et dans ce cas $\varphi_2$ est une symétrie.

\subsection{Cas général}
Notons $M_{i_1\cdots i_p}^{j_1 \cdots j_p}$ avec $i_1<i_2 < \cdots <i_p$ et $j_1<j_2 < \cdots <j_p$ le déterminant de la matrice d'ordre $p$ obtenue en ne conservant dans la matrice $a$ que les lignes $i_1,i_2 ,\cdots,i_p$ et les colonnes $j_1,j_2 ,\cdots,j_p$.
Si $i\overline{_1,i_2 ,\cdots,i_p}$ désigne le complémentaire ordonné de $i_1,i_2 ,\cdots,i_p$ dans $\left[\left[1,\cdots,n\right]\right]$, alors
\begin{propo}
$$\det (a)=\sum (-1)^{j_1+j_2+\cdots+j_p+p-1}M_{1\cdots p}^{j_1 \cdots j_p}M_{\overline{1\cdots p}}^{\overline{j_1 \cdots j_p}}$$
\end{propo}
On notera par $m_p(a)$ la matrice des cofacteurs
$$m_p(a)=(M_{i_1\cdots i_p}^{j_1 \cdots j_p})$$
les coefficients étant rangés suivant l'ordre lexicographique des indices $i_1,\cdots,i_p$ pour les lignes et $j_1,\cdots,j_p$ pour les colonnes. C'est une matrice carrée d'ordre $N=\binom{n}{p}$ inversible si et seulement si $a$ l'est et dont le déterminant est égal à $(\det(a))^{p}$.  On définit comme précédemment l'application $\varphi_p$ définie par $\varphi_p(a)=m_p(a)$. Elle définit un homomorphisme injectif de groupe
$$\varphi_p: GL(n,\K) \rightarrow GL(N,\K)$$

\end{document}